\documentclass[11pt]{article}
\usepackage{amsfonts,latexsym,rawfonts,amsmath,amssymb,amsthm,graphicx}
\numberwithin{equation}{section}
\newtheorem{theorem}{Theorem}[section]
\newtheorem{lem}[theorem]{Lemma}
\newtheorem{thm}[theorem]{Theorem}
\newtheorem{pro}[theorem]{Proposition}
\newtheorem{cor}[theorem]{Corollary}

\newtheorem{rem}[theorem]{Remark}

\def\endproof{$\hfill\Box$\\}

\title{Remarks on the nonexistence of biharmonic maps}
\author{Yong Luo}
\date{}
\begin{document}
\maketitle
\begin{abstract}
In this short note we study nonexistence result of biharmonic maps from a complete Riemannian manifold into a Riemannian manifold with nonpositive sectional curvature. Assume that $\phi:(M,g)\to (N, h)$ is a biharmonic map, where $(M, g)$ is a complete Riemannian manifold and $(N,h)$ a Riemannian manifold with nonpositive sectional curvature, we will prove that $\phi$ is a harmonic map if one of the following conditions holds:
\\(i) $|d\phi|$ is bounded in $L^q(M)$ and $ \int_M|\tau(\phi)|^pdv_g<\infty, $
for some $1\leq q\leq\infty$, $1< p<\infty$;
\\or\\(ii) $Vol(M)=\infty$ and $ \int_M|\tau(\phi)|^pdv_g<\infty, $
for some $1< p<\infty$.

In addition if $N$ has negative sectional curvature, we assume that $rank\phi(q)\geq2$ for some $q\in M$ and $\int_M|\tau(\phi)|^pdv_g<\infty, $
 for some $1< p<\infty$. These results improve  the related theorems due to Baird et al.(cf.~\cite{BFO}), Nakauchi et al.(cf.~\cite{NUG}), Maeta(cf.~\cite{Ma}) and Luo(cf.~\cite{Luo}).

\end{abstract}
\section{Introduction}
Let $(M,g)$ be a Riemannian manifold and $(N,h)$ a Riemannian manifold without boundary. For a $W^{1,2}(M,N)$ map $\phi$, the energy density of $\phi$ is defined by
$$ e(\phi)=|\nabla \phi|^2=\rm{Tr_g}(\phi^\ast h),$$
where $\phi^\ast h$ is the pullback of the metric tensor $h$. The energy functional of the mapping $\phi$ is defined as $$E(\phi)=\frac{1}{2}\int_Me(\phi)dv_g.$$  The Euler-Lagrange equation of $E$ is $\tau(\phi)=\rm Tr_g\bar{\nabla} d\phi=0$ and $\tau(\phi)$ is called the tension field of $\phi$. A map is called a harmonic map if $\tau(\phi)=0$. The theory of harmonic maps has many important applications in various fields of differential geometry, including minimal surface theory, complex geometry and so on(cf.~\cite{SY}).

Much effort has been paid in the last several decades to generalize the notion of harmonic maps. In 1983, Eells and Lemaire \cite{EL}(see also \cite{ES}) proposed to consider the bienergy functional $$E_2(\phi)=\frac{1}{2}\int_M|\tau(\phi)|^2dv_g$$ of smooth maps between Riemannian manifolds. Stationary points of the bienergy functional are called biharmonic maps. We see that harmonic maps are biharmonic maps and even more, minimizers of the bienergy functional. In 1986, Jiang\cite{Ji} derived the first and second variational formulas of the bienergy functional and studied biharmonic maps. The Euler-Lagrange equation of $E_2$ is
$$\tau_2(\phi):=-\Delta^\phi\tau(\phi)-\sum_{i=1}^mR^N(\tau(\phi), d\phi(e_i))d\phi(e_i)=0,$$
where $\Delta^\phi:=\sum_{i=1}^m(\bar{\nabla}_{e_i}\bar{\nabla}_{e_i}-\bar{\nabla}_{\nabla_{e_i}e_i})$, $\nabla$ is the Levi-Civita connection on $(M,g)$ and $\bar{\nabla}$ is the induced connection on the pull back bundle $\phi^{-1}TN$, and $R^N$ is the Riemannian curvature tensor on $N$.

The first nonexistence result of biharmonic maps was obtained by Jiang\cite{Ji}. He proved that biharmonic maps from a compact, orientable Riemannian manifold into a Riemannian manifold of nonpositive curvature are harmonic. Jiang's theorem is a direct application of the Weitzenb\"ock formula. If $\phi$ is biharmonic, then
\begin{eqnarray*}
-\frac{1}{2}\Delta|\tau(\phi)|^2&=&\langle-\Delta^\phi\tau(\phi), \tau(\phi)\rangle-|d\tau(\phi)|^2
\\&=&Tr_g\langle R^N(\tau(\phi), d\phi)d\phi, \tau(\phi)\rangle-|d\tau(\phi)|^2
\\&\leq&0.
\end{eqnarray*}
The maximum principle implies that $d\tau(\phi)=0$ and so by
$$div\langle d\phi, \tau(\phi)\rangle=|\tau(\phi)|^2+\langle d\phi, d\tau(\phi)\rangle,$$
we deduce that $div\langle d\phi, \tau(\phi)\rangle=|\tau(\phi)|^2$, then integration by parts, we have $\tau(\phi)=0$.

If $M$ is noncompact, the maximum principle is no longer applicable. In this case we can use the integration by parts argument, by choosing proper test functions. Based on this idea, Baird et al.(cf.~\cite{BFO}) proved that biharmonic maps from a complete Riemannian manifold with nonnegative Ricci curvature into a nonpositively curved manifold with finite bienergy are harmonic. It is natural to ask whether we can abandon the curvature restriction on the domain manifold and weaken the integrable condition on the bienergy. In this direction, Nakauchi et al.(cf.~\cite{NUG}) proved that biharmonic maps from a complete manifold to a nonpositively curved manifold are harmonic if($p=2$)
\\(i) $\int_M|d\phi|^2dv_g<\infty$ and $\int_M|\tau(\phi)|^pdv_g<\infty$, or
\\(ii) $Vol(M, g)=\infty$ and $\int_M|\tau(\phi)|^pdv_g<\infty.$

Later Maeta(cf.~\cite{Ma}) generalized this result by assuming that $p\geq2$. In this paper, we will further generalize this result to the following:
\begin{theorem}\label{main thm}
Let $\phi: (M, g)\to (N, h)$ be a biharmonic map from a complete Riemannian manifold $(M, g)$ into a Reimannian manifold $(N, h)$ of nonpositive sectional curvatures and  $1\leq q\leq\infty$, $1<p<\infty$.

(i) If $|d\phi|$ is bounded in $L^q(M)$ and $$\int_M|\tau(\phi)|^pdv_g<\infty,$$
then $\phi$ is harmonic.

(ii) If $Vol(M, g)=\infty$ and $$\int_M|\tau(\phi)|^pdv_g<\infty,$$
then $\phi$ is harmonic.
\end{theorem}
\begin{rem}
For a better understanding of theorem \ref{main thm} the readers could consult the papers \cite{BFO1} and \cite{LO} for examples of proper biharmonic maps(that biharmonic maps which are not harmonic).
\end{rem}
We must point out that part (ii) of theorem \ref{main thm} is in fact implictly contained in the proof of Theorem 3.1 of \cite{BFO}, which is obviously not realized by the authors of \cite{NUG} and \cite{Ma}. That is also the motivation for us to search further in this direction.

When the target manifold has negative sectional curvatures, we have
\begin{thm}\label{main2}
Let $\phi: (M, g)\to (N, h)$ be a biharmonic map from a complete Riemannian manifold $(M, g)$ into a Reimannian manifold $(N, h)$ of negative sectional curvatures and $\int_M|\tau(\phi)|^pdv_g<\infty$ for some $1<p<\infty$. Assume that there is some point $q\in M$ such that $rank\phi(q)\geq2$, then $\phi$ is a harmonic map.
\end{thm}
This theorem was proved by Oniciuc(cf.~\cite{On}) under the assumption of $|\tau(\phi)|$ is a constant.
\begin{rem}
Theorem \ref{main2} is a generalization of Theorem 1.3 in \cite{Luo}. But the right statement of Theorem 1.3 in \cite{Luo} should be added an additional assumption of $rank\phi(q)\geq2$ at some point $q\in M$.
\end{rem}

The rest of this paper is organized as follows: In section 2 we give some preliminaries on harmonic maps and biharmonic maps. In section 3 our theorems are proved. In section 4 we give some applications of our results to biharmonic submersions.
\section{Preliminaries}
\subsection{Harmonic maps and biharmonic maps}
In this section we give more details on the definitions of harmonic maps and biharmonic maps.

Let $\phi:(M, g)\to (N, h)$ be a map from an m-dimensional Riemannian manifold $(M, g)$ to an n-dimensional Riemannian manifold $(N, h)$. The energy of $\phi$ is defined by
$$E(\phi):=\frac{1}{2}\int_M|d\phi|^2dv_g.$$
The E-L equation of $E$ is
$$\tau(\phi)=\sum_{i=1}^m\{\bar{\nabla}_{e_i}d\phi(e_i)-d\phi(\nabla_{e_i}e_i)\}=0,$$
where we denote by $\nabla$ the Levi-Civita connection on $(M, g)$ and by $\bar{\nabla}$ the induced Levi-Civita connection on $\phi^{-1}TN$. $\tau(\phi)$ is called the tension field of $\phi$.

A map $\phi: (M, g)\to (N, h)$ is called a harmonic map if $\tau(\phi)=0$.

To generalize the notion of harmonic maps, in 1983 Eells and Lemaire\cite{EL}(see also \cite{ES}) proposed to consider the bienergy functional
$$E_2(\phi):=\frac{1}{2}\int_M|\tau(\phi)|^2dv_g.$$
In 1986, Jiang\cite{Ji} calculated the first and second variational formulas of the bienergy functional. The E-L equation of $E_2$ is
\begin{eqnarray}\label{equ2}
\tau_2(\phi)=-\Delta^\phi \tau(\phi)-\sum_{i=1}^mR^N(\tau(\phi), d\phi(e_i))d\phi(e_i)=0,
\end{eqnarray}
where $\{e_i, i=1,...,m\}$ is a local orthogonal frame on $M$ and $R^N$ is the Riemann curvature tensor of $(N, h)$. $\tau_2(\phi)$ is called the bitension field of $\phi$.

A map $\phi:(M, g)\to (N, h)$  is called a biharmonic map if $\tau_2(\phi)=0$.
\subsection{Gaffney's theorem}
In the subsequent section we will use the following Gaffney's theorem(cf.~\cite{Ga}).
\begin{thm}
Let $(M, g)$ be a complete Riemannian manifold. If a $C^1$ 1-form $\omega$ satisfies
that $\int_M|\omega|dv_g<\infty$ and $\int_M\delta\omega dv_g<\infty$, or equivalently, a $C^1$ vector field $X$ defined by
$\omega(Y) = \langle X, Y \rangle, (\forall Y \in TM)$ satisfies that $\int_M|X|dv_g<\infty$ and $\int_Mdiv Xdv_g<\infty$, then $$\int_M\delta\omega dv_g=\int_Mdiv Xdv_g=0.$$
\end{thm}

\section{Biharmonic maps into nonpositively curved manifolds}
In this section we will prove theorem \ref{main thm} and theorem \ref{main2}. First let's prove a lemma.
\begin{lem}\label{lem}
Assume that $\phi:(M,g)\to (N,h)$ is a biharmonic map from a complete manifold $(M,g)$ to a nonpositively curved manifold $(N,h)$ and $\int_M|\tau(\phi)|^pdv_g<\infty$ for some $p>1$. Then $|\tau(\phi)|$ is a constant and moreover $\bar{\nabla}\tau(\phi)=0$.
\end{lem}
\proof Here most part of the proof is the same with that of theorem 3.1 in \cite{BFO}. For the completeness of this paper and the convenience of the readers we give all the details here.

 Let $\epsilon>0$ and a direct computation shows that
\begin{eqnarray}\label{ine1}
&&\Delta(|\tau(\phi)|^2+\epsilon)^\frac{1}{2}\nonumber
\\&=&(|\tau(\phi)|^2+\epsilon)^{-\frac{3}{2}}(\frac{1}{2}(|\tau(\phi)|^2+\epsilon)\Delta|\tau(\phi)|^2-\frac{1}{4}|\nabla|\tau(\phi)|^2|^2).
\end{eqnarray}
Since $\nabla|\tau(\phi)|^2=2h(\bar{\nabla}\tau(\phi),\tau(\phi))$ we get
$$|\nabla|\tau(\phi)|^2|^2\leq 4(|\tau(\phi)|^2+\epsilon)|\bar{\nabla}\tau(\phi)|^2.$$
Therefore we obtain
\begin{eqnarray}
&&\frac{1}{2}(|\tau(\phi)|^2+\epsilon)\Delta|\tau(\phi)|^2-\frac{1}{4}|\nabla|\tau(\phi)|^2|^2 \nonumber
\\&\geq& \frac{1}{2}(|\tau(\phi)|^2+\epsilon)(\Delta|\tau(\phi)|^2-2|\bar{\nabla}\tau(\phi)|^2).
\end{eqnarray}
Since $\phi$ is biharmonic, from the biharmonic equation we see that
\begin{eqnarray}\label{ine3}
\frac{1}{2}\Delta|\tau(\phi)|^2=|\bar{\nabla}\tau(\phi)|^2-Tr_gR^N(\tau(\phi),d\phi,d\phi,\tau(\phi))\geq |\bar{\nabla}\tau(\phi)|^2,
\end{eqnarray}
where we used the assumption that $R^N\leq0$. Combining the above two inequalities we obtain
\begin{eqnarray}\label{ine2}
\frac{1}{2}(|\tau(\phi)|^2+\epsilon)\Delta|\tau(\phi)|^2-\frac{1}{4}|\nabla|\tau(\phi)|^2|^2\geq0.
\end{eqnarray}
From inequalities (\ref{ine1}) and (\ref{ine2}) we deduce that
$$\Delta(|\tau(\phi)|^2+\epsilon)^\frac{1}{2}\geq0.$$
Now let $\epsilon\to0$ we have that
$$\Delta|\tau(\phi)|\geq 0,$$
that $|\tau(\phi)|$ is a positive subharmonic function on $M$. By Yau's(cf.~\cite{Y}) celebrated Liouville theorem for nonnegative subharmonic functions on a complete manifold we see that if $\int_M|\tau(\phi)|^pdv_g<\infty$ for some $p>1$, then $|\tau(\phi)|$ is a constant on $M$. Moreover by (\ref{ine3}) we see that $\bar{\nabla}\tau(\phi)=0$. This completes the proof of the lemma. \endproof

 \proof of theorem \ref{main thm}. From the above lemma we see that $|\tau(\phi)|=c$ is a constant. Hence if $Vol(M)=\infty$, we must have $c=0$, this proves (ii) of theorem \ref{main thm}. To prove (i) of theorem \ref{main thm}, we distinguish two cases. If $c=0$, we are done. If not, we see that $Vol(M)<\infty$ and we will get a contradiction in the following. Define a l-form on $M$ by
$$\omega(X):=\langle d\phi(X),\tau(\phi)\rangle,~(X\in TM).$$
Then we have
\begin{eqnarray*}
\int_M|\omega|dv_g&=&\int_M(\sum_{i=1}^m|\omega(e_i)|^2)^\frac{1}{2}dv_g
\\&\leq&\int_M|\tau(\phi)||d\phi|dv_g
\\&\leq&c Vol(M)^{1-\frac{1}{q}}(\int_M|d\phi|^qdv_g)^\frac{1}{q}
\\&<&\infty ,
\end{eqnarray*}
where if $q=\infty$ we denote $\|d\phi\|_{L^\infty(M)}=(\int_M|d\phi|^qdv_g)^\frac{1}{q}$.

In addition, we consider $-\delta\omega=\sum_{i=1}^m(\nabla_{e_i}\omega)(e_i)$:
\begin{eqnarray*}
-\delta\omega&=&\sum_{i=1}^m\nabla_{e_i}(\omega(e_i))-\omega(\nabla_{e_i}e_i)
\\&=&\sum_{i=1}^m\{\langle\bar{\nabla}_{e_i}d\phi(e_i),\tau(\phi)\rangle
-\langle d\phi(\nabla_{e_i}e_i),\tau(\phi)\rangle\}
\\&=&\sum_{i=1}^m\langle \bar{\nabla}_{e_i}d\phi(e_i)-d\phi(\nabla_{e_i}e_i),\tau(\phi)\rangle
\\&=&|\tau(\phi)|^2,
\end{eqnarray*}
where in the second equality we used $\bar{\nabla}\tau(\phi)=0$. Now by Gaffney's theorem we have that
$$0=\int_M-\delta\omega=\int_M|\tau(\phi)|^2dv_g=c^2Vol(M),$$
which implies that $c=0$, a contradiction. Therefore we must have $c=0$, i.e. $\phi$ is a harmonic map. This completes the proof of theorem \ref{main thm}.
\endproof

\proof of theorem \ref{main2}. By lemma \ref{lem}, $|\tau(\phi)|=c$ is a constant. We only need to prove that $c=0$. Assume that $c\neq 0$, we will get a contradiction. Then by the biharmonic equation and the Weitzenb\"ock formula we have at $q\in M$:
\begin{eqnarray*}
0&=&-\frac{1}{2}\Delta|\tau(\phi)|^2
\\&=&-\langle\Delta^\phi\tau(\phi), \tau(\phi)\rangle-|d\tau(\phi)|^2
\\&=&\sum_{i=1}^m\langle R^N(\tau(\phi), d\phi(\partial_{x_i}))d\phi(\partial_{x_i}), \tau(\phi)\rangle-|d\tau(\phi)|^2
\\&=&\sum_{i=1}^m\langle R^N(\tau(\phi), d\phi(\partial_{x_i}))d\phi(\partial_{x_i}), \tau(\phi)\rangle,
\end{eqnarray*}
where in the first and fourth equalities we used lemma \ref{lem}. Since the sectional curvatures of $N$ is negative, we must have that $d\phi(\partial_{x_i})//\tau(\phi)$ at $q\in M$ $\forall i$, i.e. rank$\phi(q)\leq1$. Therefore we must $c=0$, a contradiction. This completes the proof of theorem \ref{main2}.
\endproof

\section{Biharmonic submersions into nonpositively curved manifolds}
In this section we give some applications of our result to biharmonic submersions.

First we recall some definitions(cf. \cite{BW}).

Assume that $\phi: (M, g)\to (N, h)$ is a smooth map between Riemannian manifolds and $x\in M$. Then $\phi$ is called \textbf{horizontally weakly conformal} if either

(i) $d\phi_x=0$, or

(ii) $d\phi_x$ maps the horizontal space $\rm \mathcal{H}_x=\{Kerd\phi_x\}^\bot$ conformally onto $T_{\phi(x)}N$, such that
$$h(d\phi_x(X), d\phi_x(Y))=\Lambda g(X, Y), (X, Y\in \mathcal{H}_x).$$

A map $\phi$ is called \textbf{horizontally weakly conformal} on $M$ if it is horizontally weakly conformal at every point of $M$. If furthermore, $\phi$ has no critical points, then we call it a \textbf{horizontally conformal submersion}. Note that if $\phi: (M, g)\to (N, h)$ is a horizontally weakly conformal map and $\rm dim M<\rm dim N$, then $\phi$ is a constant map.

If for every harmonic function $f: V\to \mathbb{R}$ defined on an open subset $V$ of $N$ with $\phi^{-1}(V)$ non-empty, the composition $f\circ\phi$ is harmonic on $\phi^{-1}(V)$, then $\phi$ is called a \textbf{harmonic morphism}. Harmonic morphisms are characterized as follows(cf. \cite{Fu, Is}).
\begin{thm}[\cite{Fu, Is}]\label{thm4}
A smooth map $\phi: (M, g)\to (N, h)$ between Riemannian manifolds is a harmonic morphism if and only if $\phi$ is both harmonic and horizontally weakly conformal.
\end{thm}
Let $\phi: (M, g)\to (N, h)$ be a submersion, then each tangent space $T_xM$ can be decomposed as follows.
\begin{eqnarray}\label{eq4}
T_xM=\mathcal{V}_x\oplus\mathcal{H}_x,
\end{eqnarray}
where $\mathcal{V}_x=\rm Ker(d\phi_x)$ is the vertical space and $\mathcal{H}_x$ is the horizontal space. If there exists a smooth function $\lambda$ on $M$ such that for each $x\in M$,
$$h(d\phi_x(X), d\phi_x(Y))=\lambda^2(x)g(X, Y),$$
then $\lambda$ is called the dilation.

When $\phi:(M^m, g)\to (N^n, h)(m>n\geq2)$ is a horizontally conformal submersion, the tension field is given by
\begin{eqnarray}\label{eq5}
\tau(\phi)=\frac{n-2}{2}\lambda^2d\phi(grad_\mathcal{H}(\frac{1}{\lambda^2}))
-(m-n)d\phi(\hat{H}),
\end{eqnarray}
where $grad_\mathcal{H}(\frac{1}{\lambda^2})$ is the horizontal component of $\rm grad(\frac{1}{\lambda^2})$ according to (\ref{eq4}), and $\hat{H}$ is the trace of the second fundamental form of each fiber given by
$$\hat{H}=\frac{1}{m-n}\sum_{i=1}^m\mathcal{H}(\nabla_{e_i}e_i),$$
where $\{e_i, i=1,...,m\}$ is a local orthonormal frame field on $M$ such that $\{e_{i}, i=1,...,n\}$ belongs to $\mathcal{H}_x$ and $\{e_{j}, j=n+1,...,m \}$ belongs to $\mathcal{V}_x$ at each point $x\in M$.

Nakauchi et al.(cf.~\cite{NUG}) and Maeta(cf.~\cite{Ma}) applied their nonexistence result of biharmonic maps to get conditions such that biharmonic submersions are harmonic morphisms. Here we generalize their result by using theorem \ref{main thm}. We have
\begin{pro}
Let $\phi:(M^m, g)\to (N^n, h)(m>n\geq2)$ be a biharmonic horizontally conformal submersion from a complete Riemannian manifold $(M, g)$ into a Riemannian manifold $(N, h)$ with nonpositive sectional curvatures and $p$ a real constant satisfying $1<p<\infty$. If
 $$ \int_M\lambda^p|\frac{n-2}{2}\lambda^2grad_\mathcal{H}(\frac{1}{\lambda^2})
-(m-n)\hat{H}|_g^pdv_g<\infty,$$
and either $\lambda$ is bounded in $L^q(1\leq q\leq\infty)$ or $Vol(M)=\infty$, then $\phi$ is a harmonic morphism.
\end{pro}
\proof By (\ref{eq5}),
$$\int_M|\tau(\phi)|_h^pdv_g=\int_M\lambda^p|\frac{n-2}{2}\lambda^2grad_\mathcal{H}(\frac{1}
{\lambda^2})
-(m-n)\hat{H}|_g^pdv_g<\infty,$$
and since $|d\phi(x)|^2=n\lambda^2(x)$, we get that $\phi$ is harmonic by theorem\ref{main thm}. Since $\phi$ is also a horizontally conformal submersion, by theorem \ref{thm4} $\phi$ is a harmonic morphism. \endproof

In particular if $\rm dimN=2$, we have
\begin{cor}
Let $\phi:(M^m, g)\to (N^n, h)(m>2)$ be a biharmonic horizontally conformal submersion from a complete Riemnnian manifold $(M, g)$ into a Rieman surface $(N, h)$ with nonpositive Gauss curvature and $p$ a real constant satisfying $1<p<\infty$. If
 $$\int_M\lambda^p|\hat{H}|_g^pdv_g<\infty,$$
and either $\lambda$ is bounded in $L^q(1\leq q\leq\infty)$ or $Vol(M)=\infty$, then $\phi$ is a harmonic morphism.
\end{cor}
Similarly, by theorem \ref{main2} we have
\begin{pro}
Let $\phi:(M^m, g)\to (N^n, h)(m>n\geq2)$ be a biharmonic horizontally conformal submersion from a complete Riemannian manifold $(M, g)$ into a Riemannian manifold $(N, h)$ with negative sectional curvatures and $p$ a real constant satisfying $1<p<\infty$. If
 $$ \int_M\lambda^p|\frac{n-2}{2}\lambda^2grad_\mathcal{H}(\frac{1}{\lambda^2})
-(m-n)\hat{H}|_g^pdv_g<\infty,$$
then $\phi$ is a harmonic morphism.
\end{pro}
This proposition generalizes (i) of proposition 4.2 in \cite{Luo}.
\quad\\

\textbf{Acknowledgement.} The author is partially supported by the Postdoctoral Science Foundation of China(No.2015M570660), and the Project-sponsored by SRF for ROCS, SEM.

{}
\vspace{1cm}\sc

Yong Luo

School of mathematics and statistics,

Wuhan university, Wuhan 430072, China

and

Max-planck institut f\"ur mathematik

In den naturwissenschaft

Inselstr.22, D-04103, Leipzig, Germany

{\tt yongluo@whu.edu.cn}~{\em or}~{\tt yongluo@mis.mpg.de}

\vspace{1cm}\sc

\end{document}